\documentclass[11pt,english]{article}
\usepackage[T1]{fontenc}
\usepackage[latin9]{inputenc}
\usepackage[a4paper]{geometry}
\usepackage{verbatim}
\usepackage{setspace}
\usepackage{amssymb}
\providecommand{\tabularnewline}{\\}

\newcommand{\citet}[1]{\cite{#1}}

\newtheorem{theorem}{Theorem}
\newtheorem{lemma}{Lemma}

\newcommand{\ds}{\displaystyle}

\newcommand{\esli}{\mbox{ \ if \ }}

\usepackage{babel}

\begin{document}
\global\long\def\nN{\mathcal{N}_{N}}
\global\long\def\ioik{\left(i_{1},\ldots,i_{k}\right)}
\global\long\def\Ls{L_{\mbox{s}}\mbox{}}
\global\long\def\Lsii{L_{\mbox{\scriptsize s},\ioik}}
\global\long\def\RR{\mathbb{R}}
\global\long\def\sE{\mathsf{E}}
\global\long\def\al{\alpha}
\global\long\def\de{\delta}
\global\long\def\mux{\boldmath{\mu}}
\global\long\def\erx{R}
\global\long\def\Lf{L_{0}\, }
\global\long\def\cdv{\varkappa}
\global\long\def\om{\omega}
\global\long\def\gn{\gamma_{N}}
\global\long\def\mud{\mu}
\global\long\def\dsp{d}
\global\long\def\sdv{b_{2}}
\global\long\def\srM{s}
\global\long\def\ZZ{\mathbb{Z}}
 \global\long\def\som{\omega'}
\global\long\def\vom{\tilde{\omega}}
 \global\long\def\cNN{\mathcal{N}_{N}}
 \global\long\def\tm{\tau_{t}^{*}}
\global\long\def\ka{k_{N}}
 \global\long\def\sP{\mathsf{P}}

\title{Time scales in large systems of Brownian particles\\ with stochastic synchronization}

\author{Anatoly MANITA%
\thanks{This work is supported by Russian Foundation of Basic Research (grant
09-01-00761). \protect\\{\it  Address\/}: Department of Probability, Faculty of Mathematics
and Mechanics, Moscow State University, 119992, Moscow, Russia. \quad
{\it E-mail\/}:~manita@mech.math.msu.su%
}}

\date{December 14, 2010}
\maketitle
\begin{abstract}
We consider a system $x(t)=(x_{1}(t),\ldots,x_{N}(t))$ consisting
of $N$ Brownian particles with synchronizing interaction between
them occurring at random time moments $\{\tau_{n}\}_{n=1}^{\infty}$.
Under assumption that the free Brownian motions and the sequence $\{\tau_{n}\}_{n=1}^{\infty}$
are independent we study asymptotic properties of the system when
both the dimension~$N$ and the time~$t$ go to infinity. We find
three time scales $t=t(N)$ of qualitatively different behavior of
the system. 
\end{abstract}
%\thispagestyle{empty}
%
\begin{comment}
\emph{Corresponding author}: 

Anatoly MANITA 

Department of Probability,

Faculty of Mathematics and Mechanics,

Moscow State University, 

119992 Moscow, Russia

~

\emph{Email}: \texttt{manita@mech.math.msu.su}

~
\end{comment}

\section{Introduction}

Mathematical models with stochastic synchronization between components
take their origin from paper~\cite{MitraMitrani} where some two-dimensional
system related with parallel computations was studied. A very good
explanation of the role of synchronizations in asynchronous parallel
and distributed algorithms can be found in~\cite{Bert-Tsits}. It~is
rather natural that further mathematical interest to such models was
moved to considerations of high dimensions and to studies of a long
time behavior. It was discovered soon \cite{man-Shch,mal-man,Malyshkin}
that it is very convenient to interpret synchronization models as
particle systems with very special interaction. It is worth to note
that in the {}``traditional'' mathematical theory of interacting
particle systems such interactions were never considered before that
time. In~\cite{man-obo} one can find a short overview of the subject. 

The present paper is a small contribution to the following general
problem: how to describe a qualitative behavior of a mutidimensional
Markov (or semi-Markov) process $x(t)=\left(x_{1}(t),x_{2}(t),\ldots,x_{N}(t)\right)$
for large $N$ and $t=t(N)$. We chose as an object of our study the
system of $N$ Brownian particles perturbed by synchronizing jumps
at some random time moments. Reasons of such choice are the following.
Synchronization models driven by Brownian motion were not studied
yet, all papers mentioned above considered random walks on lattices
or deterministic motions as non-perturbed dynamics. The second reason
is that, as it will be shown here, a Markovian synchronization model
based on the Brownian motions admits an explicit solution. This feature
give us possibility to write very short and clear proofs of our main
result on the existence of three different time stages of qualitative
behavior of the particle system. We believe that such results hold
also for very general multidimensional synchronization models. There
are already many particular examples justifying this belief. Thus
the existence of the three time stages in the long time behavior was
already proved for system with two types of deterministic particles
and pairwise stochastic synchronizations~\cite{mal-man-TVP}, for
discrete time random walks with a 3-particle anysotropic interaction~\cite{Malyshkin},
for continuous time random walks with symmetric $k$-particle synchronizations~\cite{manita-TVP-large-ident}. 

The explanatory goals of this paper force us to chose the following
organization of sections. In Section~2 we define and study a sequence
of Markov models with pairwise synchronization between particles and
constant coefficients in the front of the free dynamics and the interaction.
This lets us avoid cumbersome notation in proofs (Section~3). Section
4 is devoted to generalizations of the model of Section~2. The first
generalization to the case of coefficients varying with~$N$ is quite
straightforward and is based on a careful analysis of the proofs in
Section~3. The next extension of the main results is done for general
symmetric $k$-particle synchronizing interaction. In Subsection~4.3
we discuss generalizations to the case when epochs of synchronization
form a general renewal process and hence the particle system is no
more a Markov process. Corresponding results are obtained by using
the Laplace transform and are presented in Theorem~\ref{t-R-renewal}. 

%
\begin{comment}
This work is supported by Russian Foundation of Basic Research (grant
09-01-00761).
\end{comment}
{}

\section{Model with pairwise interaction}

\subsection{Definition and assumptions}

\label{sub:def-Br-pair}

We study a multi-dimensional stochastic process \[
x(t)=\left(x_{1}(t),x_{2}(t),\ldots,x_{N}(t)\right)\in\RR^{N},\quad\quad t\in\RR_{+},\]
which can be regarded mathematically as a special class of interacting
particle systems. But from the view point of possible applications
it would be better to consider this process as a multi-component stochastic
system.

Here $N$ is the number of particles and $x_{i}(t)\in\RR^{1}$ is
a coordinate of the $i$-th particle at time~$t$. Denote $\cNN=\left\{ 1,\ldots,N\right\} $.
To give a precise construction of the process $\left(x(t),\, t\geq0\right)$
we fix on some probability space $\left(\tilde{\Omega},\tilde{\mathcal{F}},\tilde{\sP}\right)$
\begin{description}
\item [{~}] (a) $B(t)=\left(B_{1}(t),\ldots,B_{N}(t)\right)$ --- the
$N$-dimensional standard Brownian motion,
\item [{~}] (b) a random sequence $\left\{ \tau_{n}\right\} _{n=1}^{\infty}$
of time moments \[
0=\tau_{0}<\tau_{1}<\tau_{2}<\cdots\]

\item [{~}] (c) a random initial configuration of particles $x(0)=\left(x_{1}(0),x_{2}(0),\ldots,x_{N}(0)\right)$.
\end{description}
\emph{Main assumption} is that $\left(B(t),\, t\geq0\right)$, $\left\{ \tau_{n}\right\} _{n=1}^{\infty}$
and $x(0)$ are \emph{independent}. 

We consider also another probability space $\left(\Omega',\mathcal{F}',\mathsf{P}'\right)$
corresponding to the independent sequence \begin{equation}
(i_{1},j_{1}),\,(i_{2},j_{2}),\,\ldots,\,(i_{n},j_{n}),\,\ldots\label{eq:pair-inter-0}\end{equation}
of \emph{equiprobable} ordered pairs $(i,j)$ such that $i,j\in\cNN$,
$i\not=j$. In the next we will put simply $\som=\left((i_{1},j_{1}),\,(i_{2},j_{2}),\,\ldots,\,(i_{n},j_{n}),\,\ldots\right)$
and will use coordinate functions $i_{n}(\som)=i_{n}$ and $j_{n}(\som)=j_{n}$.

Let us introduce the new probability space $\left(\Omega,\mathcal{F},\mathsf{P}\right)=\left(\tilde{\Omega}\times\Omega',\tilde{\mathcal{F}}\times\mathcal{F}',\tilde{\mathsf{P}}\times\mathsf{P}'\right)$.
\emph{By~formal definition} the process $\left(x(t),\, t\geq0\right)$
has right-continuous trajectories $\left(x(t,\om),\, t\geq0\right)$,
$\om=(\vom,\som)$, satisfying to the following conditions: \begin{eqnarray*}
x_{k}(s,\om)-x_{k}(\tau_{n}(\vom),\om) & = & \sigma\cdot\left(B_{k}(s,\vom)-B_{k}(\tau_{n}(\vom),\vom)\right),\\
 &  & \quad\quad\forall s\in[\tau_{n}(\vom),\tau_{n+1}(\vom)),\quad\forall k\in\cNN\,,\\
x_{j_{n}(\som)}(\tau_{n}(\vom),\om) & = & x_{i_{n}(\som)}(\tau_{n}(\vom)-0,\om),\quad\quad\\
x_{m}(\tau_{n}(\vom),\om) & = & x_{m}(\tau_{n}(\vom)-0,\om)\quad\quad\forall m\in\cNN\backslash\left\{ j_{n}(\som)\right\} \,.\end{eqnarray*}

The scalar parameter $\sigma>0$ is a diffusion coefficient.

\emph{Informally} speaking the dynamics of the process $x(t)$ consists
of two parts: free motion and pairwise interaction between particles.
Namely, the\emph{ interaction} is possible only at random time moments
\[
0<\tau_{1}<\tau_{2}<\cdots\]
and has the form of \emph{synchronizing jumps}: at time $\tau_{n}$
with probability $\frac{1}{N(N-1)}$ a pair of particles $(i,j)$
is chosen and the particle~$j$ jumps to the particle~$i$: \begin{equation}
\left(x_{i},x_{j}\right)\rightarrow\left(x_{i},x_{i}\right).\label{eq:pair-inter}\end{equation}
 Inside the intervals $(\tau_{k},\tau_{k+1})$ particles of the process
$x(t)$ move as independent Brownian motions with diffusion coefficient
$\sigma$ (\emph{free dynamics}).

In some sense the dynamics of the interacting particle system $x(t)$
can be considered as a perturbation of the stochastic dynamics $B(t)$.
We are interested in the question how the synchronizing interaction
will imply on a long time behavior of $x(t)$. We consider the following
limiting situations:
\begin{description}
\item [{~}] (i) $N$ is fixed, $t\rightarrow\infty$;
\item [{~}] (ii) $N\rightarrow\infty$ is fixed, $t=t(N)\rightarrow\infty$
with different choices of the time scales $t(N)$.
\end{description}
We shall mainly be concerned here with the situation (ii) which is
more important and more interesting. 

To make our considerations more transparent in all subsequent sections
we have the next assumption.
\begin{description}
\item [{Assumption~M.}] The moments $\left\{ \tau_{n}\right\} _{n=1}^{\infty}$
are epochs of a Poisson flow of intensity $\delta$, i.e., the sequence
$\left\{ \tau_{n}-\tau_{n-1}\right\} _{n=1}^{\infty}$ consists of
independent random variables, having exponential distributions: $\mathsf{P}\left(\tau_{n}-\tau_{n-1}>s\right)=\exp(-\delta s)$.
\end{description}
Assumption M implies immediately that $\left(x(t),\,\, t\geq0\right)$
is a Markov process on $\RR^{N}$ with symbolic generator\[
\sigma L_{0}^{B}+\delta\Ls,\qquad\sigma>0,\quad\delta>0\,,\]
where $L_{0}^{B}$ is a generator of the standard $N$-dimensional
Brownian motion and $\Ls$ corresponds to synchronizing jumps.

This assumption is not crucial for the validity of our asymptotic
results. In Subsection~\ref{sub:renewal} we shall discuss the case
of a general \emph{renewal process}.

\subsection{Long time behavior for fixed $N$}

We use notation $\mathcal{L}\left(\xi\right)$ for a distribution
law of a random element $\xi$. Then $\left(\mathcal{L}\left(x(t)\right),\, t\geq0\right)$
is a family of probability measures on $\left(\RR^{N},\mathcal{B}\left(\RR^{N}\right)\right)$. 

\begin{theorem} \label{t-x(t)-no-lim}

$\mathcal{L}\left(x(t)\right)$ has no limit as $t\rightarrow\infty$.

\end{theorem} 

We recall the well known fact that the Brownian motion $B(t)$ also
has no limit on distribution as $t\rightarrow\infty$. But a long
time behavior of the interacting particle system $x(t)$ strongly
differs from the behavior of $B(t)$. Indeed, let us consider an {}``improved''
process $x^{\circ}(t)$,\[
x_{i}^{\circ}(t)=x_{i}(t)-M(x(t)),\]
where $M(x):=\frac{1}{N}\sum\limits _{m=1}^{N}x_{m}\,$ is the center
of mass of the particle configuration $x=\left(x_{1},\ldots,x_{N}\right)$.
In other words, $x^{\circ}(t)$ is the particle system $x(t)$ viewed
by an observer placed in the center of mass $M(x(t))$.

\begin{theorem} \label{t-x0-lim} For any $\sigma>0,\,\,\delta>0$
the Markov process $x^{\circ}(t)$ is ergodic. Hence there exists
a probability distribution $\mu^{N}$ on $\left(\RR^{N},\mathcal{B}\left(\RR^{N}\right)\right)$
such that $\mathcal{L}\left(x^{\circ}(t)\right)\rightarrow\mu^{N}$
as $t\rightarrow\infty$.

\end{theorem} 

The idea of the proof is to show that $x^{\circ}(t)$ satisfies the
Doeblin property. Similar arguments were used in~\cite{mal-man,manita-umn}.
So we omit here the proofs of Theorems~\ref{t-x(t)-no-lim} and~\ref{t-x0-lim}.

The result of Theorem~\ref{t-x0-lim} is close to the \emph{shift-compactness}
property of measure-valued stochastic processes~\cite{dorog}.

It would be interesting to answer the following main questions. What
is a typical {}``size'' of the configuration $(x_{1},\ldots,x_{N})$
under the distribution~$\mu^{N}$? How large (with respect to $N$
is a domain where $\mu_{N}$ is supported with probability close to~1?
To do this we let the dimension~$N$ and the time $t$ grow to infinity
in order to find on which time scale $t=t(N)$ the process $x(t)$
will approach~$\mu^{N}$.

\subsection{Time scales}

In collective behavior of a particle system with synchronization we
observe a superposition of \emph{two opposite tendencies}: with the
course of time the free dynamics \emph{increases} the \emph{spread}
of the particle system while the synchronizing interaction tries to
\emph{decrease} it.

To formalize the notion of a {}``size'' or a {}``spread'' we consider
the following function on the state space

\[
V:\,\RR^{N}\rightarrow\RR_{+},\quad\quad V(x):=\frac{1}{N-1}\sum_{m=1}^{N}\left(x_{m}-M(x)\right)^{2},\]
where $M(x)$ is the center of mass as defined above. In statistics
the function $V$ is known as the empirical variance. We introduce
also the function $R_{N}:\,\RR_{+}\rightarrow\RR_{+}$ depending on
the time $t\geq0$ as \begin{equation}
R_{N}(t):=\sE\, V(x(t)).\label{eq:R-N-t}\end{equation}
It appears that the function $R_{N}(t(N))$ has completely different
asymptotic behavior for different choices of the time scale $t=t(N)$.
Before proving this result we start from the following explicit formula.

\begin{theorem} \label{t-R-expl-f} There exist a number $\varkappa>0$
such that \[
R_{N}(t)=\sigma^{2}\delta^{-1}\, l_{N}\,\left(1-\exp\left(-\delta t/l_{N}\right)\right)+\exp\left(-\delta t/l_{N}\right)R_{N}(0),\]
where $l_{N}=N(N-1)/\varkappa$. \end{theorem} This statement shows
that the function $R_{N}(t)$ satisfies to a very simple differential
equation\[
\frac{d}{dt}R_{N}(t)=\sigma^{2}-\frac{\delta}{l_{N}}\, R_{N}(t)\,.\]
So the choice of $R_{N}$ in~(\ref{eq:R-N-t}) was really good from
the point of view of subsequent asymptotic analysis.

In the next theorem we assume that $N\rightarrow\infty,\quad t=t(N)\rightarrow\infty.$

\begin{theorem}[On three time scales]\label{t-R} Let $\ds\sup_{N}\erx_{N}(0)<\infty$.
Then

\begin{tabular}{rlcl}
\textbf{I.} & If $\quad\ds\frac{t(N)}{N^{2}}\rightarrow0$, & ~then & $\quad R_{N}(t(N))\sim\,\sigma^{2}\, t(N)$.\tabularnewline
 &  &  & \tabularnewline
\textbf{II.} & If $\quad t(N)=cN^{2}/\left(\varkappa\delta\right)$, & $c>0$,~~~ & \tabularnewline
 &  & ~then & $\quad\ds R_{N}(t(N))\sim\,\ds\frac{1-e^{-c}}{c}\,\,\,\sigma^{2}\, t(N)$.\tabularnewline
 &  &  & \tabularnewline
\textbf{III.} & If $\quad\ds\frac{t(N)}{N^{2}}\rightarrow\infty$, & ~then & $\quad\ds R_{N}(t(N))\sim\left(\frac{\sigma^{2}}{\varkappa\delta}\right)\, N^{2}$.\tabularnewline
\end{tabular}

\end{theorem} 

\textbf{\emph{Remark}} 1. In case $\delta=0$ when there is no synchronization
and $x(t)$ behaves as the Brownian motion $\sigma B(t)$ the function
$R_{N}$ can be calculated explicitly: $R_{N}(t)=\sigma^{2}t$~.

\textbf{\emph{Remark}} 2. The function $f$ in the item II is strictly
decreasing: \[
f(c)=\frac{1-e^{-c}}{c}\,,\quad f(0)=1,\quad f'(c)<0,\quad f(+\infty)=0.\]

\textbf{\emph{Remark}} 3. For the pairwise synchronization~(\ref{eq:pair-inter-0})
and~(\ref{eq:pair-inter}) considered in the present section $\varkappa=2$.
Details will be given at the end of Subsection~\ref{sub:k-part-synchr}.

\subsection{Discussion of collective behavior}

We can easily observe from Theorem~\ref{t-R} that for the slowest
time scale (case I) asymptotic behavior of $R_{N}(t(N))$ is the same
as for non-perturbed dynamics. This means that a cumulative effect
of synchronization jumps on time intervals of the form $(0,o(N^{2}))$
is negligible with respect to the influence of the free dynamics.
Next observation is that on the fastest time scale (case~III) asymptotics
of $R_{N}(t(N))$ does not depend on the rate of grow of $t=t(N)$.
We interpret this phenomenon as follows: synchronization dominates
heavily on the free motion and the asymptotics $(\sigma^{2}/2\delta)N^{2}$
corresponds to the averaging of the function $f(x)$ with respect
to the limiting distribution~$\mu^{N}$. The asymptotics on the middle
time scale (case II, time intervals of the form $(c_{1}N^{2},c_{2}N^{2})$)
{}``continuously joins'' the asymptotics of the slowest and the
fastest time stages. 

As in~\cite{mal-man-TVP} one can call these consecutive stages correspondingly: 
\begin{description}
\item [{\phantom{II}I}] initial desynchronization 
\item [{\phantom{I}II}] critical slowdown of desynchronization 
\item [{III}] final stabilization. 
\end{description}

\section{Proofs }

\label{sec:proofs-2}

\subsection{Proof of Theorem~\ref{t-R-expl-f}}

Let $\Pi_{t}=\max\left\{ m:\,\tau_{m}\leq t\right\} $ and $\tm=\tau_{\Pi_{t}}$.
Obviously, $\tm=\max\left\{ \tau_{i}:\,\,\tau_{i}\leq t\right\} $.
To get $R_{N}(t)$ we shall calculate the chain of conditional expectations
as follows \[
\sE\left(\cdot\right)=\sE\,\left(\sE\left(\sE\left(\cdot\,|\,\left\{ \tau_{j}\right\} _{j=1}^{\infty}\right)\,|\,\Pi_{t}\right)\right).\]
\begin{lemma}  \label{l-usl-tau} \begin{equation}
\sE\left(V(x(t)\,|\,\,\left\{ \tau_{j}\right\} _{j=1}^{\infty}\right)=\sigma^{2}\sum_{i=0}^{\Pi_{t}-1}\ka^{\Pi_{t}-i}\cdot\left(\tau_{i+1}-\tau_{i}\right)+\sigma^{2}\cdot(t-\tm)+\ka^{\Pi_{t}}R_{N}(0)\label{eq:V-usl-tau}\end{equation}
 where $k_{N}:=\left(1-\frac{\varkappa}{N(N-1)}\right)$. \end{lemma} 

To take expectation $\sE\left(\cdot\,|\,\Pi_{t}\right)$ from the
both sides of equation~(\ref{eq:V-usl-tau}) we need to know the
joint distribution of the following form \[
\sP\left\{ \tau_{q}-\tau_{q-1}\in(x,x+dx),\,\Pi_{t}=n\right\} ,\quad\quad q\leq n\]
and the expectation of the \emph{spent waiting time} $(t-\tm)$ in
terms of~\cite{feller2}. 

\begin{lemma}  \label{l-usl-sred-interv} If Assumption M holds we
have that $\left(\Pi_{t},\, t\geq0\right)$ is the Poisson process
and \[
\sE\left(\tau_{q}-\tau_{q-1}\,|\,\Pi_{t}=n\right)=\sE\left(t-\tm\,|\,\Pi_{t}=n\right)=\frac{t}{n+1}\]
\end{lemma} 

Keeping in mind Lemmas~\ref{l-usl-tau} and~\ref{l-usl-sred-interv}
we can easily proceed with calculation of $R_{N}(t)$. Under Assumption~M
\begin{eqnarray*}
\sE\left(V(x(t)\,|\,\Pi_{t}=n\right) & = & \sigma^{2}\sum_{i=0}^{n-1}k_{N}^{n-i}\frac{t}{n+1}\,+\sigma^{2}\,\frac{t}{n+1}=\,\\
 & = & \sigma^{2}\,\frac{t}{n+1}\sum_{j=0}^{n}k_{N}^{j}\,=\,\sigma^{2}\,\frac{t}{n+1}\,\frac{1-k_{N}^{n+1}}{1-k_{N}}+\ka^{n}R_{N}(0)\,\end{eqnarray*}
Moreover, for given $t>0$ the random variable $\Pi_{t}$ has the
Poisson distribution with mean~$\delta t$. Using identity\[
\sum_{n=0}^{\infty}\frac{\alpha^{n}}{\left(n+1\right)!}=\alpha^{-1}\left(e^{\alpha}-1\right)\]
we get \begin{eqnarray*}
R_{N}(t) & = & \sum_{n=0}^{\infty}\sE\left(V(x(t)\,|\,\Pi_{t}=n\right)\,\frac{\left(\delta t\right)^{n}}{n!}\exp\left(-\delta t\right)=\\
 & = & \frac{\sigma^{2}t}{1-k_{N}}\,\left(\frac{\exp(\delta t)-1}{\delta t}-\frac{\left(\exp(k_{N}\delta t)-1\right)k_{N}}{k_{N}\delta t}\right)\exp\left(-\delta t\right)+\exp\left(-(1-k_{N})\delta t\right)R_{N}(0)=\\
 & = & \frac{\sigma^{2}}{\delta}\,\frac{1-\exp\left(-(1-k_{N})\delta t\right)}{1-k_{N}}+\exp\left(-(1-k_{N})\delta t\right)R_{N}(0)\,.\end{eqnarray*}

Putting $l_{N}=\left(1-k_{N}\right)^{-1}$ we obtain statement of
Theorem~\ref{t-R-expl-f}.

\subsection{Proof of Theorem~\ref{t-R} }

Our task is to analyze asymptotic behavior of \[
R_{N}(t(N))=\sigma^{2}\delta^{-1}\, l_{N}\,\left(1-\exp\left(-\delta t(N)/l_{N}\right)\right)+\exp\left(-\delta t(N)/l_{N}\right)R_{N}(0)\,,\qquad l_{N}=N(N-1)/\varkappa,\]
for different choices of $t=t(N)$. Let $N\rightarrow\infty$.

Case I: $t(N)/l_{N}\rightarrow0$. Then \[
R_{N}(t(N))\sim\sigma^{2}\delta^{-1}\, l_{N}\,\delta t(N)/l_{N}=\sigma^{2}t(N)\,.\]

Case II: $t(N)/l_{N}\rightarrow c\delta^{-1}$ for some $c>0$. We
have\[
R_{N}(t(N))\sim\sigma^{2}\delta^{-1}\, t(N)c^{-1}\delta\,\left(1-\exp\left(-c\right)\right)=\sigma^{2}t(N)\,\left(1-\exp\left(-c\right)\right)/c.\]

Case III: $t(N)/l_{N}\rightarrow+\infty$. Here we get\[
R_{N}(t(N))\sim\sigma^{2}\delta^{-1}\, l_{N}\sim\frac{\sigma^{2}}{\delta\varkappa}\, N^{2}\,.\]
Theorem~\ref{t-R} is proved.

\subsection{Proofs of Lemmas}

Let us introduce families of $\sigma$-algebras which are generated
$\mathcal{F}_{m}=\sigma\left(\left(x(s),\, s\leq\tau_{m}\right),\,\left\{ \tau_{i}\right\} _{i=1}^{\infty}\right)$
as follows \begin{eqnarray*}
\mathcal{F}_{m} & = & \sigma\left(\left(x(s),\, s\leq\tau_{m}\right),\,\left\{ \tau_{i}\right\} _{i=1}^{\infty}\right),\quad m=0,1,\ldots\\
\mathcal{F}_{m-} & = & \sigma\left(\left(x(s),\, s\leq\tau_{m}-0\right),\,\left\{ \tau_{i}\right\} _{i=1}^{\infty}\right),\quad m=1,2,\ldots\,.\end{eqnarray*}
Denote also $\mathcal{F}_{0-}=\sigma\left(\,\left\{ \tau_{i}\right\} _{i=1}^{\infty}\right)$.
Evidently, \[
\mathcal{F}_{0-}\subset\mathcal{F}_{0}\subset\cdots\subset\mathcal{F}_{m-}\subset\mathcal{F}_{m}\subset\mathcal{F}_{\left(m+1\right)-}\subset\mathcal{F}_{m+1}\subset\cdots\]
To prove Lemma~\ref{l-usl-tau} we shall use the following result
related with synchronizing jumps.

\begin{lemma}  \label{l-kappa} There exists $\varkappa>0$ such
that for any $m\in\mathbb{N}$\[
\sE\left(V(x(\tau_{m}))\,|\,\mathcal{F}_{m-}\right)=k_{N}\, V\left(x(\tau_{m}-0)\right),\]
where $\ds k_{N}=\left(1-\frac{\varkappa}{N(N-1)}\right)\in(0,1)$.

\end{lemma}  We postpone the proof of this lemma Subsection~\ref{sub:k-part-synchr}
where the same statement will be established for more general interactions.
Here we just note that in the case of pairwise synchronizations $\varkappa=2$.

Since the free dynamics of particles corresponds to Brownian motions
independent of the sequence $\left\{ \tau_{i}\right\} _{i=1}^{\infty}$
of synchronization moments, for any $m\in\mathbb{N}$ we have\[
\sE\left(V(x(\tau_{m+1}-0))\,|\,\mathcal{F}_{m}\right)=V(x(\tau_{m}))+\sigma^{2}\cdot(\tau_{m+1}-\tau_{m}).\]

Using Lemma~\ref{l-kappa} we get \begin{equation}
\sE\left(V(x(\tau_{m+1}))\,|\,\mathcal{F}_{m}\right)=k_{N}V(x(\tau_{m}))+k_{N}\sigma^{2}\cdot(\tau_{m+1}-\tau_{m})\,\label{eq:V-m+1-m}\end{equation}

Hence, iterating~(\ref{eq:V-m+1-m}) we come to the equation\begin{eqnarray*}
\sE\left(V(x(\tau_{m+1}))\,|\,\mathcal{F}_{m-1}\right) & = & \sE\left(\sE\left(V(x(\tau_{m+1}))\,|\,\mathcal{F}_{m}\right)\,|\,\mathcal{F}_{m-1}\right)\,=\,\\
 & = & k_{N}\sE\left(V(x(\tau_{m}))\,|\,\mathcal{F}_{m-1}\right)+k_{N}\sigma\cdot(\tau_{m+1}-\tau_{m})\end{eqnarray*}
By developing this recurrent equation we obtain \[
\sE\left(V(x(\tau_{n})\,|\,\,\left\{ \tau_{j}\right\} _{j=1}^{\infty}\right)=\sigma^{2}\sum_{i=0}^{n-1}\ka^{n-i}\left(\tau_{i+1}-\tau_{i}\right)+\ka^{n}R_{N}(0).\]
 In a similar way we get for any nonrandom $t>0$\[
\sE\left(V(x(t)\,|\,\,\left\{ \tau_{j}\right\} _{j=1}^{\infty}\right)=\sigma^{2}\sum_{i=0}^{\Pi_{t}-1}\ka^{\Pi_{t}-i}\cdot\left(\tau_{i+1}-\tau_{i}\right)+\sigma^{2}\cdot(t-\tm)+\ka^{\Pi_{t}}R_{N}(0).\]
Here we take into account that all $\tau_{j}$ have continuous distributions
and, as usually, the sign {}``='' for conditional expectations is
understood in the sense of {}``almost surely'' \cite{Shiryaev}.

This completes the proof of the Lemma~\ref{l-usl-tau}.

Lemma \ref{l-usl-sred-interv} follows from the well know facts of
renewal processes theory \cite{Gnedenko-Kovalenko,COX} or can be
verified by a direct calculation in our concrete case.

\section{Generalizations }

\subsection{Varying parameters}

Since Theorems~\ref{t-R-expl-f} and~\ref{t-R} deal with the sequence
$\left\{ x(t)=(x_{1}(t),\ldots,x_{N}(t)\right\} _{N=1}^{\infty}$
of stochastic processes is it natural to ask whenever these statements
remain true if we let the coefficients $\sigma$ and $\delta$ depend
on $N$. In other words under Assumption~M we consider a family of
Markov processes defined on the state spaces $(\RR^{N},\mathcal{B}\left(\RR^{N})\right)$
with formal generators \begin{equation}
\qquad\qquad\sigma_{N}L_{0}^{B}+\delta_{N}\Ls,\qquad\qquad\sigma_{N}>0,\quad\delta_{N}>0\,.\label{eq:s-N-d-N}\end{equation}

If we check carefully all calculations and arguments in the proofs
of Theorems~\ref{t-R-expl-f} and~\ref{t-R} we see that these proofs
are valid \emph{without any modification} for the case~(\ref{eq:s-N-d-N}).
The corresponding results are summarized in the next theorem.

\begin{theorem} \mbox{ } \label{t-R-1}
\begin{enumerate}
\item There exist a number $\varkappa>0$ (not depending on $N$) such that
\begin{equation}
R_{N}(t)=\sigma_{N}^{2}\delta_{N}^{-1}\, l_{N}\,\left(1-\exp\left(\delta_{N}t/l_{N}\right)\right)+\exp\left(-\delta_{N}t/l_{N}\right)R_{N}(0),\label{eq-rN-sN-dN-expl}\end{equation}
where $l_{N}=N(N-1)/\varkappa$. 
\item Let $N\rightarrow\infty,\quad t=t(N)\rightarrow\infty.$ Assume that
$\ds\sup_{N}\erx_{N}(0)<\infty$. There are three different time stages
in the collective behavior of the particle system: 
\end{enumerate}
\noindent \begin{center}
\begin{tabular}{c|c|c|c}
 & \textbf{I} & \textbf{II} & \textbf{III}\tabularnewline
\hline
$t(N)$ & $\alpha_{N}\rightarrow0$ & $\alpha_{N}\rightarrow c>0$ & $\alpha_{N}\rightarrow\infty$\tabularnewline
 &  &  & \tabularnewline
$R_{N}(t(N))\sim\,$ & ~$\quad\sigma_{N}^{2}t(N)\quad$~ & $\quad(1-e^{-c})c^{-1}\,\sigma_{N}^{2}\, t(N)\quad$ & $\quad\sigma_{N}^{2}\,(\varkappa\delta_{N})^{-1}N^{2}\quad$\tabularnewline
\end{tabular} 
\par\end{center}

\noindent where$\quad\ds\alpha_{N}:=\frac{\varkappa\delta_{N}\, t(N)}{N^{2}}\,\,.$

\end{theorem} 

\textbf{\emph{Remark}} 4. As it is seen from the representation~(\ref{eq-rN-sN-dN-expl})
the assumption $\ds\sup_{N}\erx_{N}(0)<\infty$ can be weakened. We
can let some growth of $\erx_{N}(0)$ in the limit $N\rightarrow\infty$
and the statement~2 of Theorem~\ref{t-R-1} still remains true.
But the conditions on the admissible growth will be different for
each time stage.

\subsection{$k$-particle synchronization}

\label{sub:k-part-synchr}

Recall our assumption~(\ref{eq:pair-inter}) on \emph{pairwise} \emph{interaction:}
we pick at random a~pair of particles $\left(x_{i},x_{j}\right)$
and move this particles as follows $(x_{i},x_{j})\rightarrow(x_{i},x_{i}).$
To study general problems of synchronization in stochastic systems
with applications to wide classes of self-organizing systems we should
face to so called \emph{multi-particle interactions. }The most general
rule of synchronizing jumps is \[
x=(x_{1},\ldots,x_{N})\rightarrow x'=(x'_{1},\ldots,x'_{N})\]
where $\left\{ x'_{1},\ldots,x'_{N}\right\} \subset\left\{ x_{1},\ldots,x_{N}\right\} ,\mbox{\ }\left\{ x'_{1},\ldots,x'_{N}\right\} \not=\left\{ x_{1},\ldots,x_{N}\right\} .$
Following the paper~\cite{manita-TVP-large-ident} we restrict ourself
here to symmetric $k$-particle interactions based on synchronizing
maps. Definition of synchronizing maps needs some preliminary notation.
First we introduce a set $\mathcal{I}:=\left\{ (i_{1},\ldots,i_{k}):\,\, i_{j}\in\nN,\,\, i_{p}\not=i_{q}\,\,(p\not=q)\,\right\} .$ 

Fix integers~$k\geq2$ and $k_{1}\geq2$, $\ldots$, $k_{l}\geq2$:
$\quad k_{1}+\cdots+k_{l}=k.$ The sequenced collection $(k_{1},\ldots,k_{l})$
will be called a \emph{signature} of interaction. Given the signature
$(k_{1},\ldots,k_{l})$ we introduce a map $\pi_{k_{1},\ldots,k_{l}}$
defined on the set~$\mathcal{I}$ as follows: $\pi_{k_{1},\ldots,k_{l}}:\,\,(i_{1},\ldots,i_{k})\mapsto\left(\Gamma_{1},\ldots,\Gamma_{l}\right),$
where $\Gamma_{j}=(g_{j},\Gamma_{j}^{\circ})$ with \begin{eqnarray*}
g_{1}=i_{1},\quad\quad\Gamma_{1}^{\circ} & = & (i_{2},\ldots,i_{k_{1}}),\\
\cdots\\
g_{l}=i_{k_{1}+\cdots+k_{l-1}+1},\quad\quad\Gamma_{l}^{\circ} & = & \left(i_{k_{1}+\cdots+k_{l-1}+1},\ldots,i_{k_{1}+\cdots+k_{l}}\right).\end{eqnarray*}
In other words the map $\pi_{k_{1},\ldots,k_{l}}$ is a special regrouping
of indices $(i_{1},\ldots,i_{k})$:

\[
(i_{1},\ldots,i_{k})=\left(i_{1},i_{2},\,\ldots,i_{k_{1}},\,\, i_{k_{1}+1},i_{k_{1}+2},\,\ldots,i_{k_{1}+k_{2}},\,\ldots,\,\ldots,i_{k}\right)\,\quad\qquad\]

\[
\begin{array}{ccccccc}
\qquad i_{1}, & \underbrace{i_{2},\ldots,i_{k_{1}},} & i_{k_{1}+1}, & \underbrace{i_{k_{1}+2},\ldots,i_{k_{1}+k_{2}}}, & \ldots & \, & \,\\
\qquad g_{1} & \Gamma_{1}^{\circ} & g_{2} & \Gamma_{2}^{\circ} & \ldots & \, & \,\end{array}\]
The map $\pi_{k_{1},\ldots,k_{l}}$\textbf{ }generates a family of
\emph{synchronizing maps} $\left\{ J_{k_{1},\ldots,k_{l}}^{\ioik},\,\ioik\in\mathcal{I}\right\} $
defined on the set $\RR^{N}$ of particle configurations:\begin{equation}
J_{k_{1},\ldots,k_{l}}^{\ioik}:\quad x=\left(x_{1},\ldots,x_{N}\right)\mapsto y=\left(y_{1},\ldots,y_{N}\right),\label{eq:sim-synchro-map}\end{equation}
where \[
y_{m}=\left\{ \begin{array}{rl}
x_{m}, & \esli m\notin\ioik,\\
x_{g_{j}}, & \esli m\in\ioik,\, m\in\Gamma_{j}\,.\end{array}\right.\]
We call the jump~(\ref{eq:sim-synchro-map}) a \emph{synchronization}
of the collection of particles $x_{i_{1}},\ldots,x_{i_{k}}$, \emph{corresponding}
to the signature $(k_{1},\ldots,k_{l})$. The configuration $J_{k_{1},\ldots,k_{l}}^{\ioik}x$
has at least $k_{1}$ particles with coordinates that are equal to
$x_{g_{1}}$, \ldots, at least $k_{l}$ particles at the point $x_{g_{l}}$.

We are ready now to define a particle system with \emph{symmetric
$k$-particle interaction} of the given signature $(k_{1},\ldots,k_{l})$.
To do this we repeat the strategy of Subsection~\ref{sub:def-Br-pair}
but with another definition of the probability space $\left(\Omega',\mathcal{F}',\sP'\right)$.
Now the space $\left(\Omega',\mathcal{F}',\sP'\right)$ corresponds
to the independent sequence \[
(i_{1}^{1},\ldots,i_{k}^{1}),\quad\ldots,\quad(i_{1}^{n},\ldots,i_{k}^{n}),\quad\ldots\]
 of \emph{equiprobable} elements of the set $\mathcal{I}$. As before
the dynamics of $x(t)$ consists of two parts: free motion and interaction.
Inside the intervals $(\tau_{k},\tau_{k+1})$ particles of the process
$x(t)$ move as independent Brownian motions with diffusion coefficient
$\sigma$ (\emph{free dynamics}). \emph{Interaction} is possible only
at the epochs $0<\tau_{1}<\tau_{2}<\cdots$ and has the following
form. At time $\tau_{n}$ with probability $\frac{1}{N(N-1)\cdots(N-k+1)}$
a set of indices $\ioik$ is chosen and the particle configuration
$\left(x_{1},\ldots,x_{N}\right)$ instantly changes to $\left(y_{1},\ldots,y_{N}\right)$
accordingly to the synchronizing map $J_{k_{1},\ldots,k_{l}}^{\ioik}$
(see~(\ref{eq:sim-synchro-map})).

Note that the pairwise interaction defined in (\ref{eq:pair-inter})
is a particular case of the symmetric $k$-particle synchronizing
interaction considered here. To see this put $k=2$, $l=1$, the signature
$(k_{1},\ldots,k_{l})=(2)$. Then $\pi_{2}:\,(i_{1},i_{2})\mapsto\Gamma_{1}$,
$\quad\Gamma_{1}=(g_{1},\Gamma_{1}^{\circ})=(i_{1},i_{2})$, $\quad g_{1}=i_{1}$,
$\Gamma_{1}^{\circ}=i_{2}$.

Let $\Ls_{,(k_{1},\ldots,k_{l})}$ denote a formal generator corresponding
to the symmetric $k$-particle interaction\emph{ }of the signature
$(k_{1},\ldots,k_{l})$. Main goal now is to generalize our results
to the Markov process $x(t)$ with generator

\[
\sigma_{N}L_{0}^{B}+\delta_{N}\Ls_{,(k_{1},\ldots,k_{l})}\,,\qquad\,\sigma_{N}>0,\quad\delta_{N}>0\,.\]
All arguments of the proof in Section~\ref{sec:proofs-2} can be
repeated as well for this case, we should only to take care about
an analog of Lemma~\ref{l-kappa}. Fortunately, the proof of Lemma~\ref{l-kappa}
for the general symmetric $k$-particle interaction\emph{ }can be
obtained by a slight modification of the proof of Lemma~2 in~\cite{manita-TVP-large-ident}.
So there is no need to repeat that proof here. We mention only the
explicit form of the constant $\varkappa$ entering in definition
of \[
k_{N}=1-\frac{\varkappa}{(N-1)N}\,.\]
It appears (see~\cite{manita-TVP-large-ident}) that $\varkappa=\sum_{j=1}^{l}k_{j}^{2}-k$.
It is easy to check that $\varkappa>0$ for any $k_{1}\geq2$, $\ldots$,
$k_{l}\geq2$ such that $k_{1}+\cdots+k_{l}=k.$

So our final conclusion is the following one.

\emph{The both statements of Theorem~\ref{t-R-1} remains true for
the particle system with symmetric $k$-particle interaction. Moreover,
$\varkappa=\varkappa(k_{1},\ldots,k_{l})=\sum_{j=1}^{l}k_{j}^{2}-k$.
The choice of the sequence $\left\{ \alpha_{N}\right\} $ is the same:
$\,\,\ds\alpha_{N}=\frac{\varkappa\delta_{N}\, t(N)}{N^{2}}\,\,.$}

Let us remark also that for the pairwise synchronization $\varkappa((2))=2.$

\subsection{Nonmarkovian model: general renewal epochs for synchronization}

\label{sub:renewal}

The next step in generalization of the model is to consider more general
sequences $\left\{ \tau_{n}\right\} $. We replace Assumption~M by
the following one.
\begin{description}
\item [{Assumption~T$_{N}$.}] For each fixed $N\in\mathcal{\mathbb{N}}$
the moments $\left\{ \tau_{n}^{(N)}\right\} _{n=1}^{\infty}$ are
epochs of some \emph{renewal} process, i.e., the sequence $\left\{ \tau_{n}^{(N)}-\tau_{n-1}^{(N)}\right\} _{n=1}^{\infty}$
consists of independent random variables, having common continuous
distribution function $F_{N}(s)=\sP\left\{ \tau_{n}^{(N)}-\tau_{n-1}^{(N)}\leq s\right\} $
satisfying $F_{N}(s)=0$ for $s\leq0$. Intervals $\tau_{n}^{(N)}-\tau_{n-1}^{(N)}$
have finite mean $\mu_{N}>0$ and variance~$d_{N}$.
\end{description}
Expected result is the following one. Consider a stochastic process
$x^{(N)}(t)=(x_{1}(t),\ldots,x_{N}(t))$, $t\geq0$, corresponding
to $N$ Brownian particles with diffusion coefficient $\sigma_{N}>0$.
Particles of $x^{(N)}(t)$ interact at epochs $\left\{ \tau_{n}^{(N)}\right\} _{n=1}^{\infty}$
according to the \emph{symmetric $k$-particle interaction} of the
signature $(k_{1},\ldots,k_{l})$. Let Assumption~T$_{N}$ holds.
\medskip{}

\noindent \textbf{Conjecture. }\emph{Assume that $\ds\sup_{N}\erx_{N}(0)<\infty$.
Let $N\rightarrow\infty,\quad t=t(N)\rightarrow\infty.$ There are
three different time stages in the collective behavior of the particle
system: }

\noindent \begin{center}
\emph{}\begin{tabular}{c|c|c|c}
 & \textbf{\emph{I}} & \textbf{\emph{II}} & \textbf{\emph{III}}\tabularnewline
\hline
\emph{$t(N)$} & \emph{$\alpha_{N}\rightarrow0$} & \emph{$\alpha_{N}\rightarrow c>0$} & \emph{$\alpha_{N}\rightarrow\infty$}\tabularnewline
 &  &  & \tabularnewline
\emph{$R_{N}(t(N))\sim\,$} & \emph{~$\quad\sigma_{N}^{2}t(N)\quad$~} & \emph{$\quad(1-e^{-c})c^{-1}\,\sigma_{N}^{2}\, t(N)\quad$} & \emph{$\quad\varkappa^{-1}\sigma_{N}^{2}\mu_{N}\, N^{2}\quad$}\tabularnewline
\end{tabular}\emph{ }
\par\end{center}

\noindent \emph{where$\quad\ds\alpha_{N}:=\frac{\varkappa\, t(N)}{\mu_{N}N^{2}}\,\,$,
}$\ds\quad\varkappa=\sum_{j=1}^{l}k_{j}^{2}-k$.\medskip{}

Evidently, $x^{(N)}(t)$ is not a Markov process. Of course, we can
not expect to have here an explicit representation for $R_{N}(t)$
as in Theorem~\ref{t-R-expl-f}. Possible proofs of the Conjecture
can be obtained by two different ways. The first one is close to Section~3
of the present paper. The idea is to represent $R_{N}(t)$ in term
of generating function of the number of renewals $\Pi_{t}$: $g(t,\zeta)=\sE\zeta^{\Pi_{t}}$
(we recall that under Assumption~T$_{N}$ $(\Pi_{t},\, t\geq0)$
is not a Poisson process). We are interested in the long time behavior
($t\rightarrow\infty$), so we take the Laplace transform of the function
$g(t,\zeta)$ in~$t$ (see~\cite[Section~3.2]{COX}), \[
g^{*}(s,\zeta)=\int_{0}^{+\infty}e^{-st}g(t,\zeta)\, dt\,,\]
to analyze its behavior for small $s$. Applying Tauberian theorems
from~\cite[Ch.~13, Section~5]{feller2} we come to the following
statement.

\begin{theorem}  \label{t-R-renewal}

\noindent Assume that $\ds\sup_{N}\erx_{N}(0)<\infty$, $N\rightarrow\infty$,
$t(N)\rightarrow\infty$. 

\emph{If }$\alpha_{N}\rightarrow0$, then ~$\quad R_{N}(t(N))\sim\sigma_{N}^{2}t(N)\, L_{1}(t(N))\,.$

\emph{If }$\alpha_{N}\rightarrow\infty$, then ~$\quad R_{N}(t(N))\sim\varkappa^{-1}\sigma_{N}^{2}\mu_{N}\, N^{2}L_{2}(t(N))\,.$

\noindent Here $L_{1}$ and $L_{2}$ are some slowly varying functions,
notation $\alpha_{N}$ is the same as in Conjecture.

\end{theorem} 

These results are slightly weaker than the corresponding items of
Theorems~\ref{t-R} or~\ref{t-R-1} but this is the best we can
do by this method. We omit details.

The second possible way of proving the above conjecture is an approach
based on embedded Markov chains. It was very effective in~\cite{mal-man-TVP}
and~\cite{manita-TVP-large-ident}. We shall devote to it a separate
paper.

~
\end{document}